\UseRawInputEncoding
\documentclass[11pt]{article}
\usepackage{amsmath}
\usepackage{amsmath,amssymb,amsthm,graphicx,mathrsfs}
\usepackage[numbers,sort&compress]{natbib}
\newcommand{\beq}{\begin{equation}}
\newcommand{\enq}{\end{equation}}

\newtheorem{Theorem}{Theorem}[section]
\newtheorem{Lemma}[Theorem]{Lemma}

\newtheorem{Definition}[Theorem]{Definition}

\newtheorem{Remark}[subsection]{Remark}

\newcommand{\benu}{\begin{enumerate}}
\newcommand{\beqa}{\begin{eqnarray}}
\newcommand{\beqan}{\begin{eqnarray*}}
\newcommand{\eay}{\end{array}}
\newcommand{\edm}{\end{displaymath}}
\newcommand{\eenu}{\end{enumerate}}
\newcommand{\eeq}{\end{equation}}
\newcommand{\eeqa}{\end{eqnarray}}
\newcommand{\eeqan}{\end{eqnarray*}}

\newcommand{\br}{\begin{Remark}}
\newcommand{\er}{\end{Remark}}

\newcommand{\bqa}{\begin{eqnarray}}
\newcommand{\eqa}{\end{eqnarray}}
\newcommand{\bqw}{\begin{eqnarray*}}
\newcommand{\eqw}{\end{eqnarray*}}

\newcommand{\non}{\nonumber}
\newcommand{\bea}{\begin{array}{cc}}
\newcommand{\ena}{\end{array}}

\begin{document}
\pagenumbering{arabic} \setcounter{page}{1}
\begin{center}

{\large \bf
Strong global attractors for a three dimensional nonclassical diffusion equation with memory}\\
 \vspace{0.20in}Yuming Qin$^{1,2,\ast}$ $\cdot$ $\ $ Xiaolei Dong$^{3}$$\ $ $\cdot$ $\ $ Alain Miranville$^{4}$$\cdot$ $\ $ Ke Wang$^{1,2}$\\

 \vspace{3mm}
\end{center}
\begin{abstract}\normalsize
In this paper, we study the strong global attractors for a three dimensional nonclassical diffusion equation with memory. First,  we prove the existence and uniqueness of strong solutions for the equations  by the Galerkin method. Then we prove the existence of global attractors for the equations in $H^2(\Omega)\cap H^1_0(\Omega)\times L^2_\mu(\mathbb{R}^+;H^2(\Omega)\cap H^1_0(\Omega))$ by the condition (C).

\end{abstract}
\vspace{3mm} \hspace{4mm}{\bf Key words: }Nonclassical diffusion equation; Memory; Strong global attractors;\\
\vspace{3mm}\hspace{2mm} {\bf AMS Subject Classifications}:35K57; 35B40; 35B41;

\section{Introduction}
\let\thefootnote\relax\footnote{*corresponding author.\\
Yuming Qin\\yuming$\_$qin@hotmail.com, yuming@dhu.edu.cn\\
Xiaolei Dong\\
xld0908@163.com\\
Alain Miranville\\
Alain.Miranville@math.univ-poitiers.fr\\
Ke Wang\\
kwang@dhu.edu.cn\\
$^{1}$\ Department of  Mathematics, Donghua University,\\
Shanghai 201620, P. R. China, \\
$^{2}$ Institute for Nonlinear Science, Donghua University, Shanghai 201620, P. R. China.\\
$^{3}$ School of Mathematics and Statistics, Zhoukou Normal University,\\
Zhoukou 466001, P. R. China,\\
$^{4}$ Laboratoire de Math$\acute{e}$matiques et Applications, UMR CNRS 7348-SP2MI,\\
Universit$\acute{e}$ de Poitiers, Boulevard Marie et Pierre Curie-T$\acute{e}$l$\acute{e}$port 2 86962,\\
Chasseneuil Futuroscope Cedex, France.\\
}

In this paper, we investigate the existence of global attractors in $H^2(\Omega)\cap H^1_0(\Omega)\times L^2_\mu(\mathbb{R}^+;H^2(\Omega)\cap H^1_0(\Omega))$ for the following nonclassical diffusion equations:
\begin{eqnarray}\left\{\begin{array}{ll}
u_t-\mu\Delta u_t-\Delta u-\int_0^\infty k(s)\Delta u(t-s)ds=g(u)+f(x),\ x\in\Omega,\ t>0,\\
u|_{\partial\Omega}=0,\ t>0,\\
u|_{t=0}=u_{0}(x),\ x\in\Omega,
\end{array}\right.\label{1.1}\end{eqnarray}
where $\Omega$ is a bounded domain in $\mathbb{R}^3$ with smooth boundary $\partial\Omega$, the nonlinearity  $g(u)$ and the given external force term
$f(x)$ satisfies some assumptions later.

The equation $\eqref{1.1}_1$ originated in the general nonclassical reaction equation:

\begin{align}u_t-\mu\Delta u_t-\Delta u+g(u)=f(x),\label{1.2}\end{align}
which was named  the nonclassical diffusion equations in the case $\mu>0$.  The equation \eqref{1.2} reduces to the
classical reaction diffusion equation when $\mu=0$.

To begin with, let us review some previous works concerning the nonclassical diffusion equations with fading memory \eqref{1.1}. For the past few decades, the existence and  asymptotic behavior of the solutions to problem \eqref{1.1} have been investigated by many authors under the different assumption conditions  and satisfies Lipschitz continuous. Wang and Zhong \cite{wz} investigated the pullback attractor of the nonclassical diffusion equations with fading memory in $H^1_0(\Omega)\times L^2_\mu(\mathbb{R}; H^1_0(\Omega))$ when the nonlinearity is critical and satisfies Lipschitz continuous. Later, Wang, Yang and Zhong \cite{wyz} studied the existence and regularity of the global attractors for the nonclassical diffusion equations with fading memory when the nonlinearity is critical and satisfies Lipschitz continuous. Chen, Chen and Tang \cite{cct} considered the fractal dimension of the global attractor for the  nonclassical diffusion equations with fading memory in a periodic boundary domain by using the fractal dimension theorem given in \cite{ckp}. Inspired by \cite{cm,cm1}, Wang and Wang studied the trajectory attractor and global attractor for the nonclassical diffusion equations with fading memory. Conti et al. \cite{cmp} considered  the optimal regularity of the global attractor for a three dimensional nonclassical diffusion equation with fading memory. Recently, Conti and Marchini \cite{conmar} investigated the finite dimensional regular exponential attractors for a three dimensional nonclassical diffusion equation with fading memory.  Xie, Li and Zeng \cite{xielz} considered the  uniform attractor for the non-autonomous nonclassical diffusion equations with fading memory when the nonlinearity $g$ satisfies the arbitrary polynomial growth.  Anh et al. \cite{att} proved the existence of the global attractors to a nonclassical diffusion equation with hereditary memory and different nonlinearities, for example, polynomial type, Sobolev type, and  exponential type, the main novelty of their result is that no restriction on the upper growth of the nonlinearity. Zhang et al. \cite{zwg}  considered the existence of the global attractors for a nonclassical diffusion equation with fading
memory by the contractive function method, and the nonlinearity satisfies arbitrary polynomial growth. Conti, Dell'Oro and Pata \cite{cop} investigated the existence of the global attractor for the nonclassical diffusion equations with hereditary memory lacking instantaneous damping by using the operator decomposition method. In other words, compared to the nonclassical diffusion equation with memory,  the main feature of the lacking instantaneous damping is that the equation does not contain a term of the form $-\Delta u$.

There have been some results \cite{wlz,swz,x,sy,liuma,ab,wq,wuzhang,zm,toan,wlq,wzl,lt} for the autonomous and non-autonomous nonclassical equation \eqref{1.2} when the nonlinearity $g(u)$ satisfies Sobolev type growth conditions. Among them, Wang, Li and Zhong \cite{wlz} investigated the existence of the strong global attractors for a nonclassical diffusion equation with subcritical nonlinearity, and proved the global attractors $\mathcal{A}_{\mu}\rightarrow\mathcal{A}_{0^+}$  in the sense of Hausdorff semidistance in $H_0^1(\Omega)$ as $\mu\rightarrow0^+$. Liu and Ma \cite{liuma} studied the existence of the exponential attractors for a three dimensional nonclassical diffusion equation in $H^2(\Omega)\cap H_0^1(\Omega)$. Dong and Qin \cite{dongqin1}  investigated the existence of pullback attractors for a nonclassical diffusion equation in $H^2(\Omega)\cap H_0^1(\Omega)$. Besides, there also have been some results \cite{at,xlz,yzxz,dongqin2} for the nonclassical equation \eqref{1.2} when the nonlinearity satisfies the arbitrary polynomial growth condition.

On the one hand, Wang, Yang and Zhong declared the result  in \cite{wyz} that the strong global attractors for equations \eqref{2.3}-\eqref{2.4} when the nonlinearity $g(s)$ is subcritical and satisfies $g'(s)\leq \ell$, $s\in \mathbb{R}$, and assumed the unknown function $u$ satisfies $\int_0^\infty e^{-\rho s}|\nabla u(-s)|^2ds\leq R$, which is s strict restriction in their proof. However, the additional conditions $g'(s)\leq \ell$ and $\int_0^\infty e^{-\rho s}|\nabla u(-s)|^2ds\leq R$ are not needed in this paper.

On the other hand, from the problem \eqref{1.1}, we can see that the equation $\eqref{1.1}_1$ contains the term $-\Delta u_t$, which is different from the classical reaction diffusion equations. For example, the solutions of the classical reaction diffusion equations have some smoothing effect, in other words, if the initial data  only belong to a weaker topology space, then the solutions of the classical reaction diffusion equations will belong to a stronger topology space with higher regularity. However, for problem \eqref{1.1}, both the initial data and the solutions belong to the same space, but the solutions  have no  higher regularity because of  appearance of $-\Delta u_t$. To overcome these difficulties,  motivated  by \cite{r}, we prove the existence and uniqueness of strong solutions to problem \eqref{1.1} and use the condition (C) in \cite{mawangzhong} that has been developed recently to prove the existence of global attractors for problem \eqref{1.1} in $H^2(\Omega)\cap H^1_0(\Omega)\times L^2_\mu(\mathbb{R}^+,H^2(\Omega)\cap H_0^1(\Omega))$ in this paper.

The paper is divided as follows. In Section 2, we give some lemmas which will be used frequently and the main result in Theorem 2.8. In Section 3, we prove the existence and uniqueness of strong solutions to problem \eqref{1.1} by using the  method in \cite{r}. In Section 4, we prove the existence of  absorbing sets for problem \eqref{1.1}. In Section 5, we first prove the asymptotic behavior of solutions for problem \eqref{1.1} by the condition (C), then we complete the proof of Theorem 2.8.

Hereafter, let letter $C$ be a general positive constant, which may vary from line to line at each step.

\section{Main Result}
\setcounter{equation}{0}

In this section, we first introduce some notations on the function spaces and norms which subjects to study the existence of the strong global
attractors to problem \eqref{1.1}. Without loss of generality, we assume $\mu=1$.

As in \cite{sy}, let

$$A=-\Delta,\ \text{with\ domain}\ D(A)=H^2(\Omega)\cap H_0^1(\Omega),$$
and consider the family of Hilbert spaces $D(A^{s/2})$, $s\in\mathbb{R}$, with the standard inner products and norms, respectively,

$$(u,u)_{D(A^{s/2})}=(A^{s/2}u,A^{s/2}u)\ \text{and}\ \|u\|_{D(A^{s/2})}=|A^{s/2}u|_2=\Big(\int_\Omega|A^{s/2}u|^2dx\Big)^{1/2}.$$

Let $H=L^2(\Omega)$, $V_1=H^1_0(\Omega)$ and $V_2=D(A)=H^2(\Omega)\cap H^1_0(\Omega)$. Denote by $(\cdot,\cdot)$ and $|\cdot|_2$ the inner
product and norm of $H$, respectively. We denote by $((\cdot,\cdot))$ and $\|\cdot\|_1$ the inner
product and norm of $V_1$, respectively. Denote by $\|\cdot\|_2$ the norm of $V_2$.

To investigate problem \eqref{1.1}, we need the following assumption conditions on $g$.

(i) As in  \cite{wzl}, for the nonlinearity function $g(u)$, we assume $g\in C(\mathbb{R})$ with $g(0)=0$, satisfies
the following dissipation and growth conditions:
\begin{align}\begin{cases}
\big|g'(u)\big|\leq C(1+|u|^4),\\
\big|g(u)\big|\leq C(1+|u|^\beta)\ \text{with}\ \beta<5,\\
 \lim\sup\limits_{|u|\rightarrow\infty}\frac{g(u)}{u}<\lambda_1,
 \end{cases}\label{2.1}\end{align}
where $\lambda_1>0$ is the first eigenvalue of $-\Delta$.

As in \cite{wzhong}, let $L^2_\mu(\mathbb{R}^+;V_\alpha)$ $(\alpha=1,2)$ be the Hilbert space of $V_\alpha$-valued functions on $\mathbb{R}^+$, which subjects to the inner and norm:
$$(\phi,\psi)_{\mu,\alpha}=\int_0^\infty\mu(s)(\phi,\psi)_\alpha\ \text{and}\ \|\phi\|^2_{\mu,\alpha}=\int_0^\infty\mu(s)\|\phi\|^2_{\alpha} ds.$$

Then we define the family of Hilbert spaces

$$\mathcal{V}_{\alpha}={V_\alpha}\times L^2_\mu(\mathbb{R}^+;V_{\alpha}),\ \alpha=1,2$$
with the norm

$$\|z\|_{\mathcal{V}_{\alpha}}=\|(u,\eta^t)\|_{\mathcal{V}_{\alpha}}=\frac{1}{2}\big(\|u\|^2_{\alpha-1}+\|u\|^2_{\alpha}+\|\eta^t\|^2_{\mu,\alpha}\big).$$

As in \cite{wzl}, the generalized Poincar$\acute{e}$ inequality is as follow:

\bqa \lambda_1\|u\|^2_\alpha\leq\|u\|^2_{\alpha+1},\ \forall u\in V_{\alpha+1}.\label{2.2}\eqa

Following \cite{Dafermos}, define
$$\eta^t(x,s)=\int_0^su(x,t-r)dr,\ s\geq0,$$
then we have
$$\eta_t^t(x,s)=u(x,t)-\eta_s^t(x,s),\ s\geq0.$$

Let $\mu(s)=-k'(s)$, which satisfies $k(+\infty)=0$ and $k(0)<\infty$, then equations \eqref{1.1} can be written as the following form:

\begin{eqnarray}\left\{\begin{array}{ll}
u_t-\mu\Delta u_t-\Delta u-\int_0^\infty \mu(s)\Delta \eta^t(s)ds=g(u)+f(x),\ x\in\Omega,\ t>0,\\
\eta_t^t(x,s)=-\eta_s^t(x,s)+u(x,t),
\end{array}\right.\label{2.3}\end{eqnarray}
which subjects to the initial-boundary value conditions
\begin{eqnarray}\left\{\begin{array}{ll}
u(x,t)=0,\ \ \eta^t(x,s)=0,\ x\in\partial \Omega,\ t\geq0,\\
u(x,0)=u_0(x),\ x\in\Omega,\\
\eta^0(x,s)=\eta_0(x,s)=\int_0^su_0(x,-r)dr,\ (x,s)\in\Omega\times\mathbb{R}^+.
\end{array}\right.\label{2.4}\end{eqnarray}

As in \cite{wzhong}, suppose that the memory kernel function $\mu(s)$ satisfies the following assumptions:
\begin{eqnarray}\left\{\begin{array}{ll}
1.\ \mu(s)\in C^1(\mathbb{R}^+)\cap L^1(\mathbb{R}^+),\ \mu(s)\geq0\ \text{and}\ \mu'(s)\leq0,\ \forall s\in\mathbb{R}^+,\\
2.\ \text{there\ exist\ a\ constant}\ \delta>0 \ \text{such\ that}\ \mu'(s)+\delta\mu(s)\leq0,\ \forall s\in\mathbb{R}^+.
\end{array}\right.\label{2.5}\end{eqnarray}

Before the discussion of problem \eqref{2.3}-\eqref{2.4}, we give the well-known condition $(C)$ (see \cite{mawangzhong}) which will be devoted to the proof of the existence of the strong global attractors to problem \eqref{2.3}-\eqref{2.4}.

\begin{Definition}(\cite{mawangzhong}) A $C^0$-semigroup $\{S(t)\}_{t\geq0}$ defined on Banach space $X$ is said to be satisfying
condition $(C)$ if for any $\varepsilon>0$ and for any bounded set $B$ of $X$, there exists a
$t(B)>0$ and a finite dimensional subspace $X_1$ of $X$ such that\\
(i) $\{\|PS(t)x\|; x\in B\}$ is bounded; and\\
(ii) $\|(I-P)S(t)x\|_X\leq\varepsilon,\ x\in B$, where $P:X\rightarrow X_1$ is a bounded projector.
\end{Definition}

\begin{Lemma} (\cite{evans}) Let $u\in W^{1,p}(0,T:X)$ for some $1\leq p\leq\infty$. Then
$$u\in C([0,T];X),$$
where
\begin{eqnarray}\|u\|_{W^{1,p}(0,T;X)}:=\left\{\begin{array}{ll}
\big(\int_0^T\|u\|_X^p+\|u'\|_X^pdt\big)^{1/p},\ \ (1\leq p< +\infty);\\
ess\sup\limits_{0\leq t\leq T}(\|u\|_X+\|u'\|_X),\ \ (p=+\infty).
\end{array}\right.\non\end{eqnarray}

\end{Lemma}

\begin{Lemma} (\cite{mawangzhong}) Let $X$ be a Banach space and the semigroup $\{S(t)\}_{t\geq0}$ is norm-to-weak
continuous on $X$. Then $\{S(t)\}_{t\geq0}$ has a global attractor in $X$ provided that the
following conditions hold true:\\
(i) $\{S(t)\}_{t\geq0}$ has a family of bounded absorbing sets $B$  in $X$;\\
(ii) $\{S(t)\}_{t\geq0}$ satisfies condition $(C)$.
\end{Lemma}
\begin{Lemma} (\cite{r}) Let $X\subset\subset H\subset Y$ be Banach spaces, with $X$ reflexive. Suppose that $u_m$ is
a sequence that is uniformly bounded in $L^2((0,T); X)$, and $\frac{du_n}{dt}$ is uniformly bounded in $L^p((0,T); Y)$,
for some $p>1$. Then there is a subsequence that converges strongly in  $L^2((0,T); H)$.
\end{Lemma}
\begin{Lemma} (\cite{r}) If $u_j\rightarrow u$ in $L^p(\Omega)$, then there is a subsequence that converges pointwise
to $u$ almost everywhere in $\Omega$.
\end{Lemma}

\begin{Lemma} (\cite{r}) Let $\mathcal{O}$ be a bounded open set in $\mathbb{R}^m$, and let $g_j$ be a sequence of functions
in $L^p(\mathcal{O})$ with $$\|g_j\|_{L^p(\mathcal{O})}\leq C,\ for\ all\ j\in \mathbb{Z}^+.$$
If $g\in L^p(\mathcal{O})$ and $g_j\rightarrow g$ almost everywhere, then $g_j\rightharpoonup g$ in $L^p(\mathcal{O})$.
\end{Lemma}

Besides, we need the following Gronwall-type lemma which was summarized in \cite{clau}.
\begin{Lemma} Let $\Phi$, $r_1$ and $r_2$  be non-negative, locally summable functions on $[\tau,+\infty)$, $\tau\in \mathbb{R}$,
which satisfy, for $\varepsilon>0$ and $0<\sigma<1$,  the differential inequality
\bqa \frac{d}{dt}\Phi(t)+\varepsilon\Phi(t)\leq r_1(t)+r_2(t)\Phi(t)^{1-\sigma}\ \text{a.\ e.}\ t\in[\tau,+\infty),\non\eqa
Assume also that, for $j=1,2,$

$$m_j=\sup\limits_{r\geq\tau}\int_r^{r+1}r_j(t)dt<+\infty.$$
Then, setting $C(\nu)=\frac{e^\nu}{1-e^{-\nu}},$
\bqa \Phi(t)\leq\frac{1}{\sigma}\Phi(\tau)e^{-\varepsilon(t-\tau)}+\frac{1}{\sigma}m_1C(\varepsilon)+(m_2C(\varepsilon\sigma))^{\frac{1}{\sigma}},\non\eqa
for any $t\in[\tau,+\infty)$.
\end{Lemma}
Now, we give the main result in this paper.
\begin{Theorem}
Assume $g$  satisfies the assumptions $(i)-(ii)$, respectively. Then problem \eqref{1.1} possesses a global attractor in $H^2(\Omega)\cap H^1_0(\Omega)\times L^2_\mu(\mathbb{R}^+;H^2(\Omega)\cap H^1_0(\Omega))$.
\end{Theorem}
We prove Theorem 2.8 by the following a series of lemmas in Sections 3-5.

\section{Existence and Uniqueness of Solutions}
\setcounter{equation}{0}
In this section, we obtain the existence and uniqueness of the strong solutions to nonclassical
reaction diffusion equations with memory in $\mathcal{V}_2$ by use of the Galerkin method.

First, we give an important lemma  about the memory kernel function $\mu(s)$.
\begin{Lemma}(\cite{wzhong}) Set $I=[0,T]$  for some $T>0$. Let the memory kernel function $\mu(s)$
satisfy \eqref{2.5}, then, for any $\eta^t\in C(I;L^2_\mu(\mathbb{R}^+;{V}_r))$, $0<r<3$, there exists a constant $\gamma>0$ such that
\bqa (\eta^t,\eta_s^t)_{\mu,r}\geq \frac{\gamma}{2}\|\eta^t\|^2_{\mu,r}.\label{3.1}\eqa
\end{Lemma}
\begin{Lemma} Define $I=[0,T]$ for some $T>0$. Assume that $g$ satisfies assumption conditions $(i)$
and $f\in L^2_{loc}(\mathbb{R},V_1)$. Then for any given initial datum $z_\tau=(u_0,\eta^0)\in \mathcal{V}_2$, there exists a unique solution $z=(u,\eta^t)\in \mathcal{V}_2$
for equations \eqref{2.3}-\eqref{2.4}. Moreover, for $i=1,2$, let $z_i(0)$ denote two initial conditions, and denote by $z_i$ the corresponding solutions to equations
\eqref{2.3}-\eqref{2.4}. Then, the following estimate holds: for $0\leq t\leq T$,
\begin{align} \|u_1(t)-u_2(t)\|^2_{2}+\|\eta_1(t)-\eta_2(t)\|^2_{\mu,2}\leq C(T)\Big(\|u_1(0)-u_2(0)\|^2_{2}+\|\eta_1(0)-\eta_2(0)\|^2_{\mu,2}\Big).\label{3.2}\end{align}
\end{Lemma}
\textbf{Proof}. Let $\{\omega_m\}$ be an orthogonal basis of $L^2(\Omega)$
which consists of eigenvectors of $A=-\Delta$ and $\lambda_m$ denote the corresponding eigenvalues, $m=1,2,...$ and $0<\lambda_1\leq \lambda_2\leq \cdot\cdot\cdot\leq \lambda_j\leq\cdot\cdot\cdot\rightarrow+\infty$ as $j\rightarrow+\infty$. At the same time, we select an orthogonal basis
$\{\xi_m\}_{m=1}^\infty$ of $L^2_\mu(\mathbb{R}^+;V_2)$.

We first prove the existence of strong solutions to equations \eqref{2.3}-\eqref{2.4} by the Galerkin method. Fix $T\geq0$, given an integer $m$, denote by $P_m$
and $Q_m$ the projections on subspaces

$$Span\{\omega_1,\omega_2,...\omega_m\}\subset V_1\ \text{and}\ Span\{\xi_1,\xi_2,...\xi_m\}\subset L^2_\mu(\mathbb{R}^+;V_2),$$
respectively. By using the Picard's theory of ordinary differential equations,  we can deduce a unique solution $z_m=(u_m,\eta_m^t)$
satisfies the following approximation equations:

\begin{align}
\begin{cases}
u_{mt}-\Delta u_{mt}-\Delta u_m-\int_0^\infty \mu(s)\Delta \eta_m^t(s)ds+g(u_m)=P_mf(x),\ x\in\Omega,\ t>0,\\
\eta_{mt}^t(x,s)=-\eta_{ms}^t(x,s)+u_m(x,t),\\
u_m(0)=P_mu_0,\ \eta^0_m(s)=P_m\eta^0,
\end{cases}\label{3.3}\end{align}
where $u_m(t)=\sum^{m}_{i=1}u_{mi}\omega_i$ and $\eta_m^t(s)=\sum_{i=1}^m\eta^t_{mi}\xi_i$.

$(i)$ \textbf{Step} 1. First, taking the inner product of \eqref{3.3} with $z_m=(u_m,\eta^t_m)$ to deduce

\begin{eqnarray} &\Big(\frac{du_{m}}{dt}-\frac{d\Delta u_{m}}{dt},u_m\Big)-\Big(\frac{d\Delta\eta_m^t}{dt},\eta_m^t\Big)_\mu-(\Delta u_m,u_m)-(\Delta\eta_m^t,\eta_{ms}^t)_\mu\nonumber\\
&\qquad+(g(u_m),u_m)=(p_mf(x),u_m).\label{3.4+}\end{eqnarray}

Now, let us build a priori estimates of  the uniform  boundedness on $[0,T]$. Multiplying $\eqref{3.3}_1$ by $u_{m}$ in $H$ and using the assumptions of $g$, we derive
\begin{eqnarray} &\frac{1}{2}\frac{d}{dt}(|u_{m}(t)|_2^2+\|u_{m}(t)\|_1^2+\|\eta_m^t\|^2_{\mu,1})
+\|u_{m}(t)\|_1^2\nonumber\\
&\qquad+\frac{\gamma}{2}\|\eta_m^t\|^2_{\mu,1}+(g(u_m(t)),u_m(t))\leq(f_m,u_{m}(t)),\label{3.4}\end{eqnarray}
where $f_m=P_mf$.

From the assumption conditions \eqref{2.1}, there exist two positive constants $\nu$ and $c_1$ such that

\begin{eqnarray}(g(u_m(t)),u_m(t))\geq-(1-\nu)\|u_m(t)\|_1^2-c_1.\label{3.5+}\end{eqnarray}

Noting that $\lambda_1|u_m|_2^2\leq\|u_m\|_1^2$  and using the $Young's$ inequality, we get

\begin{equation} |(f_m, u_m(t))|\leq \frac{1}{2\lambda_1}|f_m|_2^2+\frac{\lambda_1}{2}|u_m(t)|^2_2,\label{3.5}\end{equation}

and

\begin{align}
\|u_m(t)\|_1^2-\frac{\lambda_1}{2}|u_m(t)|^2_2=&\frac{1}{2}\|u_m(t)\|_1^2+\frac{1}{2}\|u_m(t)\|_1^2-\frac{\lambda_1}{2}|u_m(t)|^2_2\nonumber\\
\geq&\frac{1}{2}\|u_m(t)\|_1^2.\label{3.6}\end{align}

Inserting \eqref{3.5+}-\eqref{3.6} into \eqref{3.4}, we deduce

\begin{eqnarray} &&\frac{d}{dt}(|u_{m}(t)|_2^2+\|u_{m}(t)\|_1^2+\|\eta_m^t\|^2_{\mu,1})+\alpha_1\Big(|u_{m}(t)|_2^2
+\|u_{m}(t)\|_1^2+\|\eta_m^t\|^2_{\mu,1}\Big)\nonumber\\
&&\leq 2c_1+\frac{1}{\lambda}|f_m|_2^2,\label{3.7}\end{eqnarray}
where $\alpha_1=\min\{\frac{\lambda_1\nu}{1+\lambda_1},\gamma\}$.

Integrating above inequality over $[0,t]$, for $0\leq t\leq T$, we obtain
\begin{eqnarray} &&(|u_{m}(t)|_2^2+\|u_{m}(t)\|_1^2+\|\eta_m^t\|^2_{\mu,1})+\alpha_1\int_0^t\Big(|u_{m}(t)|_2^2
+\|u_{m}(t)\|_1^2+\|\eta_m^t\|^2_{\mu,1}\Big)ds\nonumber\\
&&\quad\leq (|u_{m}(0)|_2^2+\|u_{m}(0)\|_1^2+\|\eta_m^0\|^2_{\mu,1})+\Big(\frac{1}{\lambda_1}|f_m|_2^2+2c_1\Big)T,\ T>t.\nonumber\\ \label{3.8}\end{eqnarray}

Thus, we can attain the sequence $\{z_m\}_{m=1}^\infty$ is uniformly $(w.r.t.\ m)$ bounded in $L^\infty(0,T;\mathcal{V}_1)\cap L^2(0,T; \mathcal{V}_1)$.

For any $v_m\in V_1$, we gain
\begin{eqnarray} &&\Big|\int_0^\infty\mu(s)(\Delta\eta_m^t(s),v_m(t))ds\Big|\nonumber\\
&&\qquad\leq\int_0^\infty\mu(s)|\nabla\eta_m^t(s)|_2|\nabla v_m(t)|_2ds\nonumber\\
&&\qquad\leq\Big(\int_0^\infty\mu(s)|\nabla\eta_m^t(s)|_2^2ds\Big)^{1/2}
\Big(\int_0^\infty\mu(s)|\nabla v_m(t)|_2^2ds\Big)^{1/2}\nonumber\\
&&\qquad\leq\| v_m(t)\|_1\Big(\int_0^\infty\mu(s)|\nabla\eta_m^t(s)|_2^2ds\Big)^{1/2}
\Big(\int_0^\infty\mu(s)ds\Big)^{1/2},\nonumber\end{eqnarray}
which gives that $\int_0^\infty\mu(s)\Delta\eta_m^t(s)ds\in H^{-1}$.

$(ii)$ \textbf{Step} 2. Multiplying $\eqref{3.3}_1$
by $Au_m$ and integrating the result over $\Omega$, we deduce

\begin{eqnarray} &\frac{1}{2}\frac{d}{dt}(\|u_m(t)\|_1^2+\|u_m(t)\|_2^2+\|\eta_m^t\|^2_{\mu,2})+\|u_m(t)\|_2^2\nonumber\\
&\qquad+\frac{\gamma}{2}\|\eta_m^t\|^2_{\mu,2}+(g(u_m(t)),A u_m(t))\leq(f_m,Au_m(t)).\label{3.9}\end{eqnarray}

By the assumption conditions \eqref{2.1}, embedding theorem $H^{6/5}\hookrightarrow L^{2\beta}$ $(\beta<5)$ and interpolation inequality , we have
\begin{eqnarray} &&(g(u_m(t)),-\Delta u_m(t))\leq\int_{\Omega} |u_m(t)|^\beta|Au_m(t)|dx
\leq |u_m(t)|_{L^{2\beta}}^{\beta}\|u_m(t)\|_2\nonumber\\
&&\qquad\leq|u_m(t)|^\beta_{H^{6/5}}\|u_m(t)\|_2
\leq\|u_m(t)\|_1^{\frac{4\beta}{5}}\|u_m(t)\|_2^{\frac{\beta}{5}}\|u_m(t)\|_2\nonumber\\
&&\qquad\leq C\|u_m(t)\|_2^{\frac{2\beta}{5}}+\frac{1}{8}\|u_m(t)\|^2_2\leq C+\frac{1}{4}\|u_m(t)\|^2_2,\label{3.10}\end{eqnarray}

Applying the Young inequality, we derive

\begin{equation} (f_m,Au_m(t))\leq |f_m|_2^2+\frac{1}{4}\|u_m(t)\|^2_2.\label{3.11}\end{equation}

Inserting \eqref{3.10}-\eqref{3.11} into \eqref{3.9}, we obtain

\begin{eqnarray} &&\frac{d}{dt}(\|u_m(t)\|_1^2+\|u_m(t)\|_2^2+\|\eta_m^t\|^2_{\mu,2})+\alpha_1(\|u_m(t)\|_1^2+\|u_m(t)\|_2^2+\|\eta_m^t\|^2_{\mu,2})\nonumber\\
&&\quad\leq C+2|f_m|_2^2.\label{3.12}\end{eqnarray}

Integrating \eqref{3.12} on $[0,t]$ and using \eqref{3.8}, we get for any $t\in[0,T]$,

\begin{eqnarray} &&\|u_m(t)\|_1^2+\|u_m(t)\|_2^2+\|\eta_m^t\|^2_{\mu,2}+\alpha_1\int_{0}^t(\|u_m(s)\|_1^2+\|u_m(s)\|_2^2+\|\eta_m^t\|^2_{\mu,2})ds\nonumber\\
&&\quad\leq \|u_m(0)\|_1^2+\|u_m(0)\|_2^2+\|\eta_m^t(0)\|^2_{\mu,2}+CT+2|f_m|_2^2T\nonumber\\
&&\quad\leq C\Big(\|u_m(0)\|_1^2+\|u_m(0)\|_2^2+\|\eta_m(0)\|^2_{\mu,2}+|f_m|_2^2T+T\Big),\label{3.13}\end{eqnarray}
thus $\{z_m\}$ is uniformly $(w.r.t.\ m)$ bounded in $L^\infty(0,T;\mathcal{V}_2)\cap L^2(0,T; \mathcal{V}_2)$.

$(iii)$ \textbf{Step} 3. Multiplying $\eqref{3.3}_1$ by $u_{mt}$ and integrating the result over $\Omega$, and using the Cauchy inequality, it follows that

\begin{equation} \frac{d}{dt}\Big(\|u_m(t)\|_1^2+2 G(u_m(t))\Big)+|u_{mt}(t)|_2^2+2\|u_{mt}(t)\|_1^2\leq |f_m|_2^2+\|\eta_m^t\|^2_{\mu,1}.\label{3.14}\end{equation}

Integrating \eqref{3.14} on $[0, t]$, we deduce

\begin{eqnarray} &&\|u_m(t)\|_1^2+\int_0^t\big(|u_{mt}(s)|_2^2+2\|u_{mt}(s)\|_1^2\big)ds\nonumber\\
&&\quad\leq \int_0^t(|f_m(s)|_2^2+\|\eta_m^t\|^2_{\mu,1})ds+\|u_m(0)\|_1^2\nonumber\\
&&\qquad+2\int_\Omega G(u_m(0))dx-2\int_\Omega G(u_m(t))dx.\label{3.15}\end{eqnarray}

Using the interpolation theorem, we have

\begin{eqnarray} \int_\Omega G(u_m(0))dx&\leq&\int_\Omega |G(u_m(0))|dx\leq C(1+|u_m(0)|_{L^\infty}^{\beta+2})\nonumber\\
&\leq& C(1+\|u_m(0)\|_1^{\frac{\beta+2}{2}}\|u_m(0)\|_2^{\frac{\beta+2}{2}}),\label{3.16}\\
\int_\Omega G(u_m(t))dx&\leq&\int_\Omega |G(u_m(t))|dx\leq C(1+|u_m(t)|_{L^\infty}^{\beta+2})\nonumber\\
&\leq& C(1+\|u_m(t)\|_1^{\frac{\beta+2}{2}}\|u_m(t)\|_2^{\frac{\beta+2}{2}}).\label{3.17}\end{eqnarray}

Inserting \eqref{3.16}-\eqref{3.17} into \eqref{3.15}, we derive
\begin{eqnarray} &&\|u_m(t)\|_1^2+\int_0^t\big(|u_{mt}(s)|_2^2+2\|u_{mt}(s)\|_1^2\big)ds\nonumber\\
&&\qquad\leq \int_0^t\|\eta_m^t\|^2_{\mu,1}ds+|f_m|_2^2T+\|u_m(0)\|_1^2+C(1+\|u_m(0)\|_2^{\beta+2}+\|u_m(t)\|_2^{\beta+2})\nonumber\\
&&\qquad\leq C(T).\label{3.18}\end{eqnarray}

Then it follows that $\int_0^T\big(|u_{mt}(t)|_2^2+\|u_{mt}(t)\|_1^2\big)dt$ is uniformly $(w.r.t\ m)$ bounded.

$(vi)$ \textbf{Step} 4. Multiplying $\eqref{3.3}_1$ by $Au_{mt}$ and integrating the result over $\Omega$, and using the Cauchy inequality, it follows that

\begin{eqnarray} &&\frac{1}{2}\frac{d}{dt}\|u_m(t)\|_2^2+\|u_{mt}(t)\|_1^2+\|u_{mt}(t)\|_2^2=\int_\Omega f_m(x)Au_{mt}(t)dx\nonumber\\
&&\quad-\int_\Omega g(u_m(t))Au_{mt}(t)dx+\int_\Omega\int_0^\infty\mu(s)\Delta\eta_m^t(s)\Delta u_{mt}(t)dsdx.\label{3.19}\end{eqnarray}

Using the Cauchy inequality, we have

\begin{eqnarray} \int_\Omega f_m(x)Au_{mt}(t)dx&\leq&|f_m|_2\|u_{mt}(t)\|_2\leq \frac{3}{2}|f_m(x)|_2^2+\frac{1}{6}\|u_{mt}(t)\|_2^2,\label{3.20}\\
\int_\Omega\int_0^\infty\mu(s)\Delta\eta_m^t(s)\Delta u_{mt}(t)dsdx&\leq& \frac{3}{2}\|\eta_m^t\|^2_{\mu,2}+\frac{1}{6}\|u_{mt}(t)\|_2^2,\label{3.21}\\
\int_\Omega g(u_m(t))Au_{mt}(t)dx&\leq&\frac{3}{2}|g(u_m(t))|_2^2+\frac{1}{6}\|u_{mt}(t)\|_2^2,\label{3.22}\\
|g(u_m(t))|_2^2&\leq&\int_\Omega|g(u_m(t))|^2dx\leq C(1+|u_m(t)|_{L^\infty}^{2\beta+2})\nonumber\\
&\leq& C(1+\|u_m(t)\|_2^{2\beta+2}).\label{3.23}\end{eqnarray}

Inserting \eqref{3.20}-\eqref{3.23} into \eqref{3.19} and integrating the result on $[0,t]$, and using the result of \emph{Step} 2, we get

\begin{eqnarray} &&\|u_m(t)\|_2^2+\int_0^t\big(2\|u_{mt}(s)\|_1^2+\|u_{mt}(s)\|_2^2ds\nonumber\\
&&\qquad\leq \|u_m(0)\|_2^2+3\int_0^t\|\eta_m^t\|^2_{\mu,2}ds+|f_m|_2^2T+C\int_0^t(1+\|u_m(s)\|_{2}^{2\beta+2})ds\nonumber\\
&&\qquad\leq C\Big(\|u_m(0)\|_1^2+\|u_m(0)\|_2^2+\|\eta_m(0)\|^2_{\mu,2}+|f_m|_2^2T+T\Big)\nonumber\\
&&\qquad\leq C(T).\label{3.24}\end{eqnarray}

Then it follows that $\int_0^T\big(\|u_{mt}(t)\|_1^2+\|u_{mt}(t)\|_2^2\big)dt$ is uniformly $(w.r.t.\ m)$ bounded.

Using  \eqref{3.13} and \eqref{3.23}, we get

\begin{equation} \int_0^T|g(u_m(s))|_2^2ds\leq C\int_0^T(1+\|u_m(s)\|_2^{2\beta+2})ds\leq C,\nonumber\end{equation}
implies
\begin{equation}\{g(u_m)\}\ \text{is\ uniformly\ (w.r.t.\ m) bounded\ in}\ L^2(0,T; L^2(\Omega)).\label{3.25}\end{equation}

Therefore, we may extract a subsequence, still relabelled by $(u_m,\eta_m^t)$, and by compactness result
(see \cite{t}) such that

\begin{eqnarray}\left\{\begin{array}{ll}
u_m\rightharpoonup u,\  \text{weakly}\star\ \text{in}\ L^\infty(0,T;V_2);\\

u_m\rightharpoonup u,\ \text{weakly}\ \text{in}\ L^2(0,T;V_2);\\

u_{mt}\rightharpoonup u_t,\ \text{weakly}\ \text{in}\ L^2(0,T;V_2);\\

\eta_m^t\rightharpoonup \eta^t,\  \text{weakly}\star\ \text{in}\ L^\infty(0,T;V_2);\\
\eta_m^t\rightharpoonup \eta^t,\  \text{weakly}\ \text{in}\ L^2(0,T;V_2);\\

g(u_m)\rightharpoonup \xi,\text{weakly}\ \text{in}\ L^2(0,T;L^2(\Omega)).
\end{array}\right.\label{3.26}\end{eqnarray}

Then, by Lemma 2.4, it follows that $u_m\rightarrow u $ in $L^2(0,T; L^2(\Omega))$. Applying Lemma 2.5, there exists a
subsequence, still relabelled by $u_m$, such that $u_m \rightarrow u $ for almost every $(x,t)\in \Omega_T=\Omega\times (0,T)$. Using the continuity of $g$, that
$g(u_m(x,t))\rightarrow g(u(x,t))$ for almost every $(x,t)\in \Omega_T$. Combining with the bound on $g(u_m(x,t))$ in $L^2(\Omega_T)$ given in \eqref{3.25}, we apply Lemma 2.6 to deduce that $g(u_m)\rightharpoonup g(u)$ in $L^2(\Omega_T)$. By the uniqueness of weak limits, it follows that $\xi=g(u)$.

Therefore, we take limit for \eqref{3.4+} and find that $z=(u,\eta^t)\in L^\infty(\tau,T; \mathcal{V}_2)$ is a solution to problem \eqref{2.3}-\eqref{2.4}, which is equivalent to
\begin{eqnarray} &&\Big(\frac{du}{dt}-\frac{d\Delta u}{dt},v\Big)-\Big(\frac{d\xi^t}{dt},\Delta\eta^t\Big)_\mu-(\Delta u,v)-(\Delta\eta^t,\xi_{s}^t)_\mu\nonumber\\
&&\qquad+(g(u),v)=(f(x),v),\label{3.28}\end{eqnarray}
for every $Z=(v,\xi^t)\in \mathcal{V}_2$ and for almost every $t\in[0,T]$.

Next, we prove continuity of solutions to problem \eqref{2.3}-\eqref{2.4} from $[0,T]$ to $\mathcal{V}_2$.

Using $u\in L^\infty(0,T;V_2)\cap L^2(0,T;V_2) $, we have \begin{equation} \int_0^T|g(u(s))|_2^2ds\leq C\Big(\int_0^T\big(|u(s)|_2^2+\|u(s)\|_2^{2\beta+2}\big)ds\Big)\leq C.\nonumber\end{equation}

Then, by Lemma 2.1, $v= u+ Au\in L^\infty(0,T;H)$ and $v'=-Au-g(u)-\int_0^\infty \mu(s)\Delta \eta^t(s)ds+f(t)\in L^2(0,T;H)$,
 which imply $u\in C^0([0,T];V_2)$.

The continuity of $\eta^t$ had been verified in \cite{zwg}, that is, $\eta^t\in C(0,T;L^2_\mu(\mathbb{R}^+; V_2)).$

To show that $u(0)=u_0$ and $\eta^0=\eta_0$, we adapt the same method employed in \cite{r}. Choosing
$\varphi\in C^1(0,T; V_1)$ with $\varphi(T)=0$ in limit equation \eqref{3.4+}. We can integrate
\eqref{3.28} over $[0,T]$ and then carry out integrating by parts on the first two terms to deduce
\begin{eqnarray}&&-\int_0^T(u,\varphi')ds-\int_0^T(\nabla u,\nabla\varphi')ds+\int_0^T(\nabla u,\nabla \varphi)ds
+\int_0^T(\nabla\eta^t,\nabla\varphi)_\mu ds\nonumber\\
&&\qquad+\int_0^T(g(u),\varphi)ds=(u(0),\varphi(0))+(\nabla u(0),\nabla\varphi(0))+\int_0^T(f(x),\varphi)ds,\label{3.29}\end{eqnarray}

If we repeat the same process, then from \eqref{3.4+}, we get
\begin{eqnarray} &&-\int_0^T(u_m,\varphi')ds-\int_0^T(\nabla u_m,\nabla\varphi')ds+\int_0^T(\nabla u_m,\nabla \varphi)ds
+\int_0^T(\nabla\eta_m^t,\nabla\varphi)_\mu ds\nonumber\\
&&\qquad+\int_0^T(g(u_m),\varphi)ds=(u_m(0),\varphi(0))+(\nabla u_m(0),\nabla\varphi(0))+\int_0^T(P_mf(x),\varphi)ds.\label{3.30}\end{eqnarray}

Taking the limit $m\rightarrow +\infty$ in \eqref{3.30}, we have
\begin{eqnarray} &&-\int_0^T(u,\varphi')ds-\int_0^T(\nabla u,\nabla\varphi')ds+\int_0^T(\nabla u,\nabla \varphi)ds
+\int_0^T(\nabla\eta^t,\nabla\varphi)_\mu ds\nonumber\\
&&\qquad+\int_0^T(g(u),\varphi)ds=(u_0,\varphi(0))+(\nabla u_0,\nabla\varphi(0))+\int_0^T(f(x),\varphi)ds,\label{3.31}\end{eqnarray}
because of $u_m(0)=u_0$. Since $\varphi(0)$ is arbitrary, a comparison of \eqref{3.29} and \eqref{3.31} gives that $u(0)=u_0$.

Assume $z_1=(u_1,\eta_1^t)$ and $z_2=(u_2,\eta_2^t)$ are two solutions to problem \eqref{2.3}-\eqref{2.4} with the initial data $z_{10}$ and $z_{20}$.

Let $u=u_1-u_2$, $\eta^t=\eta_1^t-\eta_2^t$, then $z=(u,\eta^t)$ satisfies the following equation

\begin{align}
\begin{cases}
u_t-\Delta u_t-\Delta u-\int_0^\infty \mu(s)\Delta \eta^t(s)ds+(g(u_1)-g(u_2))=0,\ x\in\Omega,\ t>0,\\
\eta_t^t(x,s)=-\eta_s^t(x,s)+u(x,t).
\end{cases}\label{3.32}\end{align}

Multiplying $\eqref{3.32}_1$ by $Au$ in $H$ and integrating by parts, we find

\begin{eqnarray}&&\frac{1}{2}\frac{d}{dt}(\|u(t)\|_1^2+\|u(t)\|_2^2+\|\eta^t\|^2_{\mu,2})
+\|u(t)\|_2^2+\frac{\gamma}{2}\|\eta^t\|^2_{\mu,2}\nonumber\\
&&\quad+(g(u_1(t))-g(u_2(t)),-\Delta u(t))\leq0.\label{3.33}\end{eqnarray}

Using the assumption conditions \eqref{2.1} of $g$ and the Young inequality, we get

\begin{eqnarray} &&\frac{d}{dt}(\|u(t)\|_1^2+\|u(t)\|_2^2+\|\eta^t\|^2_{\mu,2})\leq h(t)\Big(\|u(t)\|_1^2
+\|u(t)\|_2^2+\|\eta^t\|^2_{\mu,2}\Big),\label{3.34}\end{eqnarray}
where $h(t)=C(1+\|u_1\|_2^{2\beta+2}+\|u_2\|_2^{2\beta+2}).$

By the Gronwall inequality, it follows

\begin{eqnarray}\|u(t)\|_1^2+\|u(t)\|_2^2+\|\eta^t\|^2_{\mu,2}\leq (\|u(0)\|_1^2+\|u(0)\|_2^2+\|\eta(0)\|^2_{\mu,2})e^{\int_0^th(s)ds},
\ \forall t\geq 0.\label{3.35}\end{eqnarray}

Therefore we conclude the continuous dependence of the solutions on the initial data,
and in particular, the uniqueness when $(u(0),\eta(0))=(0,0)$. The proof is thus complete.\hfill $\Box$

Therefore, we can define the solution semigroup $\{S(t)\}_{t\geq0}$ in the space $\mathcal{V}_2$ as follows:
\begin{eqnarray} S(t): \mathcal{V}_2\rightarrow \mathcal{V}_2,\ S(t)(u_0,\eta^0)=(u(t),\eta^t),
\ \forall t\geq0.\label{3.36}\end{eqnarray}
\section{Absorbing Sets}
\setcounter{equation}{0}
In this section, we will establish the existence of global absorbing sets for a three dimensional nonclassical diffusion
equation with memory. Throughout this subsection, we always assume that the initial data belong
to a bounded set of corresponding suitable space.

Next, we show the existence of the global absorbing balls $\{\mathcal{B}(t)\}$ in $\mathcal{V}_1$ and $\mathcal{V}_2$ for problem \eqref{2.3}-\eqref{2.4}.

\begin{Lemma}Let $g$ satisfy the assumption conditions $(i)$, then the  semigroup $\{S(t)\}_{t\geq 0}$ possesses a family of global absorbing balls $\{\mathcal{B}(0,\rho_1)\}$ in $\mathcal{V}_1$ with the center zero and $\rho^2_1=\frac{2}{\alpha_1}\Big(2c_1+\frac{1}{\lambda}|f|_2^2\Big)$,
\bqa \mathcal{B}(0,\rho_1)=\Big\{(u, \eta^t)\in \mathcal{V}_1:\ \|(u, \eta^t)\|^2_{\mathcal{V}_1}\leq \rho^2_1\Big\}.\label{4.1}\eqa
\end{Lemma}

\textbf{Proof}. Multiplying \eqref{2.3} by $u$ in $H$ and using the assumptions of $g$, we derive
\begin{eqnarray} &\frac{1}{2}\frac{d}{dt}(|u(t)|_2^2+\|u(t)\|_1^2+\|\eta^t\|^2_{\mu,1})
+\|u(t)\|_1^2\nonumber\\
&\qquad+\frac{\gamma}{2}\|\eta^t\|^2_{\mu,1}+(g(u(t)),u(t))=(f,u(t)).\label{4.2}\end{eqnarray}

Applying the condition assumptions \eqref{2.1}, there exist two positive constants $\nu$ and $c_1$ such that

\begin{eqnarray}(g(u(t)),u(t))\geq-(1-\nu)\|u(t)\|_1^2-c_1.\label{4.3}\end{eqnarray}

Using the Poincare inequality $\lambda_1|u|_2^2\leq\|u\|_1^2$  and  the $Young's$ inequality, we get

\begin{equation} |(f, u(t))|\leq \frac{1}{2\lambda_1}|f|_2^2+\frac{\lambda_1}{2}|u(t)|^2_2,\label{4.4}\end{equation}
and

\begin{align}
\|u(t)\|_1^2-\frac{\lambda_1}{2}|u(t)|^2_2=&\frac{1}{2}\|u(t)\|_1^2+\frac{1}{2}\|u(t)\|_1^2-\frac{\lambda_1}{2}|u(t)|^2_2\nonumber\\
\geq&\frac{1}{2}\|u(t)\|_1^2.\label{4.5}\end{align}

Substituting \eqref{4.3}-\eqref{4.5} into \eqref{4.2}, we deduce

\begin{eqnarray} &&\frac{d}{dt}(|u(t)|_2^2+\|u(t)\|_1^2+\|\eta^t\|^2_{\mu,1})+\alpha_1\Big(|u(t)|_2^2
+\|u(t)\|_1^2+\|\eta^t\|^2_{\mu,1}\Big)\nonumber\\
&&\quad\leq 2c_1+\frac{1}{\lambda}|f|_2^2,\label{4.6}\end{eqnarray}
where $\alpha_1=\min\{\frac{\lambda_1\nu}{1+\lambda_1},\gamma\}$.

By the Gronwall inequality for \eqref{4.6}, we conclude
\begin{eqnarray} |u(t)|_2^2+\|u(t)\|_1^2+\|\eta^t\|^2_{\mu,1}\leq (|u_0|_2^2+\|u_0\|_1^2+\|\eta^0\|^2_{\mu,1})e^{-\alpha_1t}+ \frac{1}{\alpha_1}\Big(2c_1+\frac{1}{\lambda}|f|_2^2\Big).\label{4.7}\end{eqnarray}

If $|u_0|_2^2+\|u_0\|_1^2+\|\eta^0\|^2_{\mu,1}\leq M_1$, then there exists a $t_0=\frac{1}{\alpha_1}ln\frac{M_1\alpha_1}{\Big(2c_1+\frac{1}{\lambda}|f|_2^2\Big)}$ such that
\begin{eqnarray}|u_0|_2^2+\|u_0\|_1^2+\|\eta^0\|^2_{\mu,1}\leq\frac{1}{\alpha_1}\Big(2c_1+\frac{1}{\lambda}|f|_2^2\Big),\ \forall t\geq t_0.\nonumber\end{eqnarray}

Thus

\begin{eqnarray} |u(t)|_2^2+\|u(t)\|_1^2+\|\eta^t\|^2_{\mu,1}\leq \frac{2}{\alpha_1}\Big(2c_1+\frac{1}{\lambda}|f|_2^2\Big):=\rho^2_1,\ \forall t\geq t_0.\label{4.8}\end{eqnarray}
The proof is now complete. \hfill$\Box$
\begin{Lemma}Let $g$ satisfy the assumption conditions $(i)$, then the  semigroup $\{S(t)\}_{t\geq 0}$ possesses a family of global absorbing balls $\{\mathcal{B}(0,\rho_2)\}$ in $\mathcal{V}_2$ with the center zero and $\rho_2^2=C\Big(1+\frac{1}{\sigma}m_1C(\beta_2)+(m_2C(\beta_2\sigma))^{\frac{1}{\sigma}}\Big)$,
\bqa \mathcal{B}(0,\rho_2)=\Big\{(u, \eta^t)\in \mathcal{V}_2:\ \|(u, \eta^t)\|^2_{\mathcal{V}_2}\leq \rho_2^2\Big\}.\label{4.9}\eqa
\end{Lemma}
\textbf{Proof}. Multiplying \eqref{2.3} by $Au$ and integrating the result over $\Omega$, we attain

\bqa &&\frac{1}{2}\frac{d}{dt}(\|u(t)\|_1^2+\|u(t)\|_2^2+\|\eta^t\|^2_{\mu,2})
+\|u(t)\|_2^2+\frac{\gamma}{2}\|\eta^t\|^2_{\mu,2}\nonumber\\
&&\quad\leq(g(u(t)),-\Delta u(t))+(f, -\Delta u(t)).\label{4.10}\eqa

By the assumption conditions \eqref{2.1}, embedding theorem $H^{6/5}\hookrightarrow L^{2\beta}$$(\beta<5)$ and interpolation inequality, we have
\bqa &&(g(u(t)),-\Delta u(t))\leq\int_{\Omega} |u(t)|^\beta|Au|dx
\leq |u(t)|_{L^{2\beta}}^{\beta}\|u(t)\|_2\non\\
&&\quad\leq|u(t)|^\beta_{H^{6/5}}\|u(t)\|_2
\leq\|u(t)\|_1^{\frac{4\beta}{5}}\|u(t)\|_2^{\frac{\beta}{5}}\|u(t)\|_2
\leq C\|u(t)\|_2^{\frac{2\beta}{5}}+\frac{1}{4}\|u(t)\|^2_2,\label{4.11}\eqa
and

\bqa (f(x),Au(t))\leq |f(x)|_2^2+\frac{1}{4}\|u(t)\|^2_2.\label{4.12}\eqa

Inserting \eqref{4.11}-\eqref{4.12} into \eqref{4.10} and using \eqref{2.2}, we obtain

\bqa &&\frac{d}{dt}(\|u(t)\|_1^2+\|u(t)\|_2^2+\|\eta^t\|^2_{\mu,2})
+\beta_2(\|u(t)\|_1^2+\|u(t)\|_2^2+\|\eta^t\|^2_{\mu,2})\non\\
&&\quad\leq C(\|u(t)\|_1^2+\|u(t)\|^2_2+\|\eta^t\|^2_{\mu,2})^{\frac{\beta}{5}}+|f|_2^2, \label{4.13}\eqa
where $\beta_2=\min\{\frac{\lambda_1}{1+\lambda_1},\gamma\}$.

Applying Lemma 2.7, for any $t\in[t_0,+\infty)$, we get

\bqa &&\|u(t)\|_1^2+\|u(t)\|_2^2+\|\eta^t\|^2_{\mu,2}\leq\frac{1}{\sigma}(\|u(t_0)\|_1^2+\|u(t_0)\|_2^2+\|\eta^t(t_0)\|^2_{\mu,2})
e^{-\beta_2(t-t_0)}\non\\
&&\qquad+\frac{1}{\sigma}m_1C(\beta_2)+(m_2C(\beta_2\sigma))^{\frac{1}{\sigma}}\nonumber\\
&&\quad\leq\frac{1}{\sigma}(\|u(t_0)\|_1^2+\|u(t_0)\|_2^2+\|\eta^t(t_0)\|^2_{\mu,2})
e^{-\beta_2(t-t_0)}+C_2,\label{4.14}\eqa
where $\sigma=1-\frac{\beta}{5}$, $m_1=|f|_2^2$ and $m_2=C$.

Noting that $\frac{1}{\sigma}(\|u(t_0)\|_1^2+\|u(t_0)\|_2^2+\|\eta^t(t_0)\|^2_{\mu,2})\leq M_2$, ($M_2$ is a positive constant only depending on initial date), then there exists a $t_1=\max\{t_0, t_0+\frac{1}{\beta_2}ln\frac{M_2}{C_2}\}$ such that
\begin{eqnarray}M_2e^{-\beta_2(t-t_0)}\leq C_2,\ \forall t\geq t_1,\nonumber\end{eqnarray}
which along with \eqref{4.14}, impies

\begin{eqnarray} \|u(t)\|_1^2+\|u(t)\|_2^2+\|\eta^t\|^2_{\mu,2}\leq 2C_2:=\rho_2^2,\ \forall t\geq t_1.\label{4.15}\end{eqnarray}

The proof is thus complete. \hfill$\Box$

\section{Global Attractors in $\mathcal{V}_2$}
\setcounter{equation}{0}
In this section, we use the condition (C) to prove the asymptotic  compactness of semigroup.

\begin{Lemma} Under the assumptions of Theorem 2.8, the semigroup $\{S(t)\}_{t\geq0}$ satisfies condition (C).
\end{Lemma}

\textbf{Proof}. Let $\{\omega_m\}$ be an  orthogonal basis of $H_m$
which consists of eigenvectors of $A=-\Delta$. $\lambda_m$ denote the corresponding eigenvalues, $m=1,2,...$ and $0<\lambda_1\leq \lambda_2\leq \cdot\cdot\cdot\leq \lambda_j\leq\cdot\cdot\cdot\rightarrow+\infty$ as $j\rightarrow+\infty$. Then $\{\omega_m\}$ is also an orthogonal basis of
 $V_2$. We write $H_m = span\{\omega_1,\cdot\cdot\cdot,\omega_m \}$ and $P_m:V_2\rightarrow H_m$ is the orthogonal projection.

Since $f\in H$,  it can be also assumed that $m$ is large enough such that

\begin{equation}|(I-P_m)f|_{H}\leq\frac{\varepsilon}{2}\label{5.1}.\end{equation}

For any $(u, \eta^t)\in \mathcal{V}_2$, we decompose $(u,\eta^t)$ into two parts, respectively,
\begin{align}\begin{cases}
 u=P_mu+(I-P_m)u\overset{\Delta}=v+w,\\
 \eta^t=P_m\eta^t+(I-P_m)\eta^t\overset{\Delta}=\zeta^t+\xi^t,
 \end{cases}\label{5.2}\end{align}
where $$\zeta^t(x,s)=\int_0^sv(x,t-r)dr,\ \xi^t(x,s)=\int_0^sw(x,t-r)dr,\ \forall s\geq0.$$

Multiplying \eqref{2.3} by $Aw$ and integrating over $\Omega$, we conclude
\begin{eqnarray} &&\frac{1}{2}\frac{d}{dt}(|A^{\frac{1}{2}}w(t)|_2^2+|Aw(t)|_2^2+\|\xi^t(t)\|^2_{\mu,2})+|Aw(t)|_2^2+
\frac{\gamma}{2}\|\xi^t(t)\|^2_{\mu,2}\nonumber\\
&&\quad\leq-(g(u),Aw)+(f,Aw).\label{5.3}\end{eqnarray}

Using the assumption conditions \eqref{2.1} on $g$, H$\ddot{o}$lder's inequality and \eqref{5.1}, we have

\begin{eqnarray} &&\frac{1}{2}\frac{d}{dt}(|A^{\frac{1}{2}}w(t)|_2^2+|Aw(t)|_2^2+\|\xi^t(t)\|^2_{\mu,2})+|Aw(t)|_2^2+
\frac{\gamma}{2}\|\xi^t(t)\|^2_{\mu,2}\nonumber\\
&&\quad\leq (1+\|u(t)\|^4_{L^\infty})|\nabla u(t)|_{2}|A^{\frac{1}{2}}w(t)|_2+\frac{\varepsilon}{2}|Aw(t)|_2\nonumber\\
&&\quad\leq (1+|A^{\frac{1}{2}} u(t)|^2_2|A u(t)|^2_2)|A^{\frac{1}{2}} u(t)|_2|A^{\frac{1}{2}}w(t)|_2+\frac{\varepsilon}{2}|Aw(t)|_2.\nonumber\end{eqnarray}

Utilizing the Young inequality, $\lambda_{m+1}\|w(t)\|^2_{1}\leq\|w(t)\|_{2}^2$ and Lemmas 4.1-4.2, and taking $\varepsilon$ small enough, we derive

\begin{eqnarray} &&\frac{d}{dt}(\|w(t)\|_{1}^2+\|w(t)\|_{2}^2+\|\xi^t(t)\|^2_{\mu,2})
+\delta(\|w(t)\|_{1}^2+\|w(t)\|_{2}^2+\|\xi^t(t)\|^2_{\mu,2})\nonumber\\
&&\quad\leq \frac{\varepsilon^2}{2}+\frac{C}{\lambda_{m+1}},\ \forall t\geq t_1,\label{5.4}\end{eqnarray}
where $\delta=\min\{\frac{\lambda_1}{1+\lambda_1},\gamma\}.$

Using the Gronwall inequality (\cite{qin,qin1,qin2}) for \eqref{5.4} and letting $m\rightarrow+\infty$, we deduce

\begin{eqnarray}&&\|w(t)\|_{1}^2+\|w(t)\|_{2}^2+\|\xi^t(t)\|^2_{\mu,2}\nonumber\\
&&\quad\leq
(\|w(t_1)\|_{1}^2+\|w(t_1)\|_{2}^2+\|\xi^{t_1}\|^2_{\mu,2})e^{-\delta (t-t_1)}+\frac{\varepsilon^2}{\delta}\nonumber\\
&&\quad\leq \rho^2_2e^{-\delta (t-t_1)}+\frac{\varepsilon^2}{\delta},\nonumber \ \forall t\geq t_1.
\end{eqnarray}

Taking $t_2-t_1\geq\delta log\big(\frac{\rho^2_2}{\varepsilon^2}\big)$, we have

\begin{eqnarray}\|w(t)\|_{1}^2+\|w(t)\|_{2}^2+\|\xi^t(t)\|^2_{\mu,2}\leq
(1+\frac{1}{\delta})\varepsilon^2,\label{5.5} \ \forall t\geq t_2.
\end{eqnarray}

Therefore, we can conclude the semigroup $\{S(t)\}_{t\geq0}$ satisfies the condition (C). \hfill$\Box$

Combining Lemmas 3.1, 4.2 with Lemma 5.1, we finish the proof of Theorem 2.8.

\textbf{Acknowledgments:}
This paper was in part supported by the National Natural Science Foundation
of China with contract number 12171082, the fundamental research funds for
the central universities with contract numbers 2232022G-13,2232023G-13 and by a grant from science and technology commission of Shanghai municipality.


\begin{thebibliography}{lllp}

\bibitem{ab}Anh C. T., Bao T. Q.: Pullback attractors for a class of non-autonomous nonclassical diffusion equations.
Nonlinear Analysis TMA. 73, 399-412(2010)
%
%
%
\bibitem{att}Anh C. T., Thanh  D. T. P., Toan N. D.: Global attractors for nonclassical diffusion equations with hereditary memory and a new class of nonlinearities.  Ann. Polon. Math. 119, 1-21(2017)

\bibitem{at} Anh C. T., Toan  N. D.: Pullback attractors for nonclassical diffusion equations in noncylindrical domains. International Journal of Mathematics and Mathematical Sciences. (2020). DOI:10.1155/2012/875913.

\bibitem{ckp}  Celebi A. O., Kalantarov  V. K., Polat  M.: Attractors for the generalized Benjamin-Bona-Mahony
Equation. J. Differential Equations. 157, 439-451(1999)
%
\bibitem{cct} Chen T., Chen  Z., Tang Y. B.: Finite dimensionality of global attractors for a non-classical reaction diffusion equation with memory. Appl. Math. Lett. 25, 357-362(2012)

\bibitem{cm} Chepyzhov V. V., Miranville  A.: Trajectory and global attractors of dissipative hyperbolic equations with
memory. Commun. Pure Appl. Anal. 4, 115-142(2005)

\bibitem{cm1} Chepyzhov V. V., Miranville A.: On trajectory and global attractors for semilinear heat equations with fading
memory. Indiana University Math. 55, 119-167(2006)


\bibitem{clau}  Claudio G., Munoz  J. E., Pata V.:  Global attractors
for a semilinear hyperbolic equation in viscoelasticity. J. Math. Anal. Appl. 260, 83-99(2001)

\bibitem{conmar} Conti M., Marchini E. M.: A remark on nonclassical diffusion equations with memory. Appl. Math. Optim. 73, 1-21(2016)

\bibitem{cmp} Conti M., Marchini  E. M., Pata V.: Nonclassical diffusion with memory.  Math. Methods Appl. Sci. 38, 948-958(2015)

\bibitem{cop} Conti M., Dell'Oro  F. D., Pata V.: Nonclassical diffusion with memory lacking instantaneous damping.  Comm. Pure Appl. Math.
19, 2035-2050(2020)

\bibitem{Dafermos} Dafermos C. M.: Asymptotic stability in viscoelasticity. Arch. Ration. Mech. Anal. 37, 297-308(1970)


\bibitem{dongqin1} Dong  X., Qin  Y.: Strong pullback attractors for a nonclassical diffusion equation. Differ. Contin. Dyn. Syst. B. 2022. DOI:10.3934/dcdsb.2021313.



\bibitem{dongqin2} Dong  X., Qin Y.: Finite fractal dimension of pullback attractors for a nonclassical diffusion
equation. AIMS Math. 7(5), 8064-8079(2022)


\bibitem{evans} Evans L. C.: Partial Differential Equations, Vol. 19. American Mathematical Society, (1998)


\bibitem{lt} Lee J., Toi V. M.: Attractors for nonclassical diffusion equations with dynamic boundary conditions. Nonlinear Anal. 195, 111737(2020)


\bibitem{liuma}  Liu Y., Ma  Q.: Exponential attractor for a nonclassical diffusion equation. Electron, J. Differential Equations. 9, 1-7(2009)



\bibitem{mawangzhong} Ma Q., Wang  S., Zhong C.: Necessary and sufficient conditions for the existence of global attractor
for semigroup and application. Indiana University Mathematics Journal. 51(6), 1541-1559(2002)



\bibitem{qin1}  Qin Y.: Integral and Discrete Inequalities and Their Applications, Volume I: Nonlinear
Inequalities, Birkh\"{a}user, (2016)
\bibitem{qin2} Qin Y.: Integral and Discrete Inequalities and Their Applications: Volume II: Nonlinear
Inequalities, Birkh\"{a}user, (2016)
\bibitem{qin}  Qin Y.: Analytic Inequalities and Their Applications in PDEs, Birkh$\ddot{a}$user, (2017)


\bibitem{r}  Robinson J. C.: Infinite Dimentional Dynamical System. Cambridge University Press, (2001)

%
%
%
\bibitem{swz}  Sun C., Wang  S.,  Zhong C.: Global attractors for a nonclassical diffusion equation. Acta Math. Sin. Engl. Ser. 23,  1271-1280(2007)
%
\bibitem{sy} Sun  C., Yang M.: Dynamics of the nonclassical diffusion equations.  Asymptot. Anal. 59, 51-81(2008)



\bibitem{t}  Temam R.: Infinite Dimensional Dynamical Systems in Mechanics and Physics, volume
68. Appl. Math. Sci. Springer-Verlag, New York, second edition, (1997)



\bibitem{toan} Toan N. D.: Existence and long-time behavior of variational solutions to a class of nonclassical diffusion
equations in noncylindrical domains. Acta Math Vietnam. 41, 37-53(2016)



\bibitem{wlq} Wang Y., Li  P., Qin Y.: Upper semicontinuity of uniform attractors for nonclassical diffusion equations. Bound. Value Probl. 2017:84(2017)

\bibitem{wlz} Wang S., Li  D., Zhong C.: On the dynamics of a class of nonclassical parabolic equations. J. Math. Anal. Appl. 317, 565-582(2006)

\bibitem{wq} Wang Y., Qin Y.: Upper semicontinuity of pullback attractors for nonclassical diffusion equations. J. Math. Phys. 5, 022701(2010)

\bibitem{wyz} Wang X., Yang  L., Zhong C.: Attractors for the nonclassical diffusion equations with fading memory. J. Math. Anal. Appl. 362, 327-337(2010)

\bibitem{wz}  Wang Y., Zhong C.: On the existence of pullback attractors for nonautonomous reaction diffusion equations. Dyn. Syst. 23, 1-16(2008)


\bibitem{wzhong} Wang  X., Zhong C.; Attractors for the non-autonomous nonclassical diffusion equations with fading memory.
Nonlinear Anal. Theory Methods Appl. 71(11), 5733-5746(2009)

\bibitem{wzl} Wang Y., Zhu  Z., Li P.: Regularity of pullback attractors for nonautonomous nonclassical diffusion equations. J. Math. Anal. Appl. 459, 16-31(2018)

\bibitem{wuzhang} Wu  H., Zhang Z.: Asymptotic regularity for the nonclassical diffusion equation with
lower regular forcing term. Dyn. Syst. 26, 391-400(2011)

\bibitem{x} Xiao Y.:  Attractors for a nonclassical diffusion equation.  Acta Math. Appl. Sin. 18, 273-276(2002)

\bibitem{xielz} Xie Y., Li  Y., Zeng Y.: Uniform attractors for nonclassical diffusion equations with memory.  J. Funct. Spaces. http://dx.doi.org/10.1155/2016/5340489.
%
\bibitem{xlz} Xie Y., Li  J., Zhu K.: Upper semicontinuity of attractors for nonclassical diffusion equations with
arbitrary polynomial growth. Adv. Difference Equ. 2021:75(2021)


\bibitem{yzxz} Yuan J., Zhang S., Xie Y., Zhang J.: Exponential attractors for nonclassical diffusion equations with arbitrary
polynomial growth nonlinearity. AIMS Mathematics. 6(11), 11778-11795(2021)


\bibitem{zm} Zhang  Y., Ma Q.: Exponential attractors of the nonclassical diffusion equations with lower regular
forcing term. International Journal of Modern Nonlinear Theory and Application. 3, 15-22(2014)
%
\bibitem{zwg} Zhang Y., Wang  X., Gao C.: Strong global attractors for nonclassical diffusion equation with fading memory. Adv. Difference Equ.
(2017):163(2017)




\end{thebibliography}
\end{document}